\documentclass[12pt]{article}
\usepackage{amsmath,amssymb,fullpage}
\usepackage[applemac]{inputenc}

\newtheorem{definition}{Definition}[section]

\newtheorem{theorem}[definition]{Theorem}

\numberwithin{equation}{section}

\def\RB{\mathbb{R}}

\def\FC{\mathcal{F}}

\def\1B{\text{1\!\!I}}

\begin{document}
\title{A comparison theorem for  backward SPDEs with jumps}

\date{17 February 2014}

\author{
Bernt \O ksendal$^{1,2}$ \and Agn\`es Sulem$^3$, \and
Tusheng Zhang$^{4}$}
\footnotetext[1]{
Center of Mathematics for Applications (CMA), Dept. of Mathematics, University of Oslo,
P.O. Box 1053 Blindern, N--0316 Oslo, Norway,
email: {\tt oksendal@math.uio.no}. The research leading to these results has received funding from the European Research Council under the European Community's Seventh Framework Programme (FP7/2007-2013) / ERC grant agreement no [228087].}
\footnotetext[2]{
Norwegian School of Economics and Business Administration,
Helleveien
30, N--5045 Bergen, Norway.}
\footnotetext[3]{ INRIA Paris-Rocquencourt, Domaine de Voluceau, Rocquencourt, BP 105, Le Chesnay Cedex, 78153, France and 
 Universit\'e Paris-Est, F-77455 Marne-la-Vall\'ee, France
and CMA,University of Oslo,  email: {\tt agnes.sulem@inria.fr}}
\footnotetext[4]{ School of Mathematics, University of Manchester, Oxford Road, Manchester M13 9PL, United Kingdom and School of Mathematics, University of Science and Technology of China, Hefei, China, email:{\tt Tusheng.zhang@manchester.ac.uk}}

\maketitle
\begin{center}
In honor of Professor M. Fukushima on occasion of his 80th birthday
\end{center}

\begin{abstract}
In this paper, we obtain a comparison theorem for backward stochastic partial differential equation (SPDEs) with jumps.
We apply it to introduce space-dependent convex risk measures as a model for risk in large systems of interacting components.
\end{abstract}

\noindent{\bf Key Words}: Stochastic partial differential equations (SPDEs), comparison theorem, backward stochastic partial differential equations with jumps, convex risk measures.
\vskip 0.3cm
\noindent \emph{MSC(2010)}\textbf{:} Primary 60H15 Secondary 93E20,
35R60.

\section{Introduction and framework}
\setcounter{equation}{0}
There are several papers dealing with comparison theorems for backward stochastic partial differential equations (BSPDEs). One of the first seems to be the paper \cite{MY}. The results of that paper were subsequently extended (still for linear BSPDEs only) in the paper \cite{DM}. Other related papers are \cite{DQT} and also our own paper \cite{OSZ} (for reflected BSPDE).

The paper which seems to be closest to ours is \cite{MYZ}. Here more general non-linear BSPDEs are considered, and a comparison theorem is proved for such equations by exploiting the relation between BSPDEs and coupled systems of  forward-backward SDEs (FBSDEs).

Our paper also deals with quite general non-linear BSPDEs, but it differs from \cite{MYZ} in several ways:

(i) First, our paper includes jumps.

(ii) Second, our BSPDEs are slightly different. They have stronger conditions on the second order term, but allow more general drift terms.

(ii) Third, our method is different, being based on an approximation technique.

\vskip 0.4cm
Let $B_t = B_t(\omega), t\geq 0$ be a Brownian motion
and $\tilde{N}(dt,dz):=N(dt, dz)-\nu (dz)dt$
an independent compensated Poisson random measure on a filtered probability space $(\Omega, \FC, \{ \FC_t\}_{0 \leq t \leq T}, P)$,  where $\nu$ is the L\'evy measure associated with the Poisson
measure $N(\cdot, \cdot)$ on $[0,\infty)\times \RB$.
Let $D$ be a bounded domain in $\RB^d$. Denote by $A(t)$ the following   second order differential operator on $D$ equipped with the Dirichlet boundary condition:
$$ A(t)=\sum\limits_{i,j=1}^{d}{\partial\over {\partial x_i}}(a_{ij}(t,x){\partial\over {\partial x_j}}),$$
 where $a=(a_{ij}(t,x)): [0,T]\times D\rightarrow  \RB^{d\times d}$  is a measurable, symmetric matrix-valued function.
Set $L=L^2(\RB, \nu)$. Let $ b(t,x,u,v,Z, r(\cdot))$ be a measurable mapping from $[0,T]\times\RB\times \RB\times \RB^d\times\RB\times L$ into $\RB$.  Let $\beta(t)=(\beta_1(t), ..., \beta_d(t)), t\geq 0$ be a given progressively measurable $\RB^d$-valued stochastic process. From now on, if $u(t,x)$ is a function of $(t,x)$, we sometimes write $u(t)$ for the function $u(t, \cdot )$.  Consider the solution of the following  backward stochastic partial differential equation (BSPDE):
\begin{eqnarray}\label{1.1}
du(t,x) &=&-A(t) u(t) dt -b(t,x,u(t,x),\nabla u(t,x), Z(t,x), r(t,x,\cdot))dt+Z(t,x)dB_{t}
\nonumber \\
&&-<\beta(t), \nabla Z(t,x)>dt+\int_\RB r(t,x, z)\tilde{N}(dt,dz), t\in (0,T)\nonumber\\
u(T,x) &=&\phi(x) \quad a.s.
\end{eqnarray}
Here $\phi$ is an $\FC_T$-measurable $H:=L^2(D)$-valued random variable. Let $V$ be the Sobolev space $H_0^{1,2}(D)$ and $V^*$ be its dual.
\begin{definition}
An adapted random field $u(t,x)$ is said to be the solution of the BSPDE (\ref{1.1}) if

(i) $u\in D([0,T];H)\cap L^2(\Omega\times [0,T]; V)$,

(ii) for $t>0$, $u$ satisfies the following equation  almost surely in $V^*$
\begin{eqnarray}\label{1.2}
u(t,x) &=&\phi(x)+\int_t^TA(s) u(s) ds +\int_t^T b(s,x,u(s,x),\nabla u(s,x), Z(s,x), r(s,x,\cdot))ds\nonumber\\
&&-\int_t^TZ(s,x)dB_{s}+\int_t^T<\beta(s), \nabla Z(s,x)>ds-\int_t^T\int_\RB r_1(s,x, z)\tilde{N}(ds,dz)\nonumber\\
&&
\end{eqnarray}
\end{definition}
\noindent{\bf Remark}.
Equation (\ref{1.2}) is equivalent to that for any $\psi\in V$,
 \begin{eqnarray}\label{1.3}
<u(t,\cdot ), \psi> &=&<\phi(\cdot), \psi>+\int_t^T<A(s)u(s), \psi> ds \nonumber\\
&& +\int_t^T <b(s,\cdot,u(s,\cdot),\nabla u(s,\cdot), Z(s,\cdot), r(s,\cdot,\cdot)),\psi> ds\nonumber\\
&&-\int_t^T<Z_1(s,\cdot),\psi> dB_{s}-\int_t^T\int_D<\beta(s), \nabla \psi(x)>Z(s,x)dxds\nonumber\\
&&-\int_t^T\int_\RB <r_1(s,\cdot, z),\psi>\tilde{N}(ds,dz)
\end{eqnarray}
almost surely.
\vskip 0.4cm

The aim of this paper is  to prove a comparison theorem for the above  BSPDEs with jumps.

\section{Main result}

 Introduce the following assumptions:

\noindent ({\bf A.1}). There exists $ \delta_1>0$ and $0<a<1$ such that
\begin{equation}\label{a1}
\sum\limits_{i,j=1}^{d}a_{ij}(t,x)z_iz_j\geq (\frac{1}{2a} |\beta|^2(t)+\delta_1) |z|^2, \,\,\forall\, z\in \RB^d\,\,\,\mbox{and} \,\,\,x\in D
\end{equation}

\noindent ({\bf A.2}). There exists $C>0$ such that
\begin{equation}\label{a2}
|b(t,x,u_1, v_1, Z_1,r)-b(t,x,u_2, v_2, Z_2,r)|\leq  C(|u_1-u_2|+|Z_1-Z_2|+|v_1-v_2|)
\end{equation}

\noindent ({\bf A.3}).
\begin{equation}\label{a3}
b(t,x,u,v,Z,r_1)-b(t,x,u,v,Z,r_2)\leq  \int_{\RB}(r_1(z)-r_2(z))\lambda(t,z)\nu(dz),
\end{equation}
where $0\leq \lambda(t,z)\leq C(1\wedge |z|)$.

\vskip 0.3cm
\noindent For $i=1,2$, consider BSPDEs:
\begin{eqnarray}\label{2.1}
du_i(t,x) &=&-A(t) u_i(t) dt -b_i(t,x,u_i(t,x),\nabla u_i(t,x), Z_i(t,x),r_i(t,x,\cdot))dt+Z_i(t,x)dB_{t}\nonumber \\
&&-<\beta(t), \nabla Z_i(t,x)>dt+\int_\RB r_i(t,x, z)\tilde{N}(dt,dz), t\in (0,T)
\nonumber \\
u_i(T,x) &=&\phi_i(x) \quad a.s.,
\end{eqnarray}
See e.g. \cite{OPZ} for  information about BSPDEs with jumps.

The following theorem is the main result of this paper.
\begin{theorem}(Comparison theorem)
Assume (A.1), (A.2) (A.3) hold for one of the coefficients $b_i$, say,  $b_2$. If  $\phi_1(x)\leq \phi_2(x)$ and
$$ b_1(t,x,u_1(t,x),\nabla u_1(t,x), Z_1(t,x),r_1(t,x,\cdot))\leq b_2(t,x,u_1(t,x),\nabla u_1(t,x),Z_1(t,x),r_1(t,x,\cdot)),$$
then we have
$u_1(t,x)\leq u_2(t,x), x\in D$, a.e. for every $t\in [0, T]$.
\end{theorem}
\textbf{Proof}.  For $n\geq 1$, define functions $\psi_n(z)$, $f_n(x)$ as follows (see \cite{MP}).
\begin{equation}\label{2.16}
\psi_n(z)=\left\{\begin{array}{lll}0& \mbox{if $z\leq 0$},\\
2nz &\mbox{if $0\leq z\leq \frac{1}{n}$},\\
2&\mbox{if $ z>\frac{1}{n}$}.\end{array}\right.
\end{equation}
\begin{equation}\label{2.17}
f_n(x)=\left\{\begin{array}{ll}0& \mbox{if $x\leq 0$},\\
\int_0^xdy\int_0^y\psi_n(z)dz &\mbox{if $x>0$}.\end{array}\right.
\end{equation}
We have
\begin{equation}\label{2.18}
f_n^{\prime}(x)=\left\{\begin{array}{lll}0& \mbox{if $x\leq 0$},\\
nx^2 &\mbox{if $0\leq x\leq \frac{1}{n}$},\\
2x-\frac{1}{n} &\mbox{if $x>\frac{1}{n}$}.
\end{array}\right.
\end{equation}
Also $f_n(x)\uparrow (x^+)^2$ as $n\rightarrow \infty$. For $h\in K:=L^2(D)$, set
$$F_n(h)=\int_Df_n(h(x))dx.$$
$F_n$ has the following derivatives for $h_1, h_2\in K$,
\begin{equation}\label{2.19}
F_n^{\prime}(h)(h_1)=\int_Df_n^{\prime}(h(x))h_1(x)dx,
\end{equation}

\begin{equation}\label{2.20}
F_n^{\prime\prime}(h)(h_1, h_2)=\int_Df_n^{\prime\prime}(h(x))h_1(x)h_2(x)dx.
\end{equation}

Applying Ito's formula we obtain
\begin{eqnarray}\label{2.21}
&&F_n(u_1(t)-u_2(t))\nonumber\\
&&=F_n(\phi_1-\phi_2)+\int_t^TF_n^{\prime}(u_1(s)-u_2(s))(A(s)(u_1(s)-u_2(s)))ds\nonumber\\
&&+\int_t^TF_n^{\prime}(u_1(s)-u_2(s))(b_1(s,u_1(s), \nabla u_1(s), Z_1(s), r_1(s))\nonumber\\
&&\quad\quad \quad -b_2(s,u_2(s), \nabla u_2(s), Z_2(s), r_2(s)))ds\nonumber\\
&&+\int_t^TF_n^{\prime}(u_1(s)-u_2(s))(<\beta(s), \nabla Z_1(s)-\nabla Z_2(s)>)ds\nonumber\\
&&-\int_t^TF_n^{\prime}(u_1(s)-u_2(s))(Z_1(s)-Z_2(s))dB_s\nonumber\\
&&-\frac{1}{2}\int_t^TF_n^{\prime\prime}(u_1(s)-u_2(s))(Z_1(s)-Z_2(s),Z_1(s)-Z_2(s) )ds\nonumber\\
&&-\int_t^T\int_\RB\bigg \{ F_n(u_1(s-)-u_2(s-)+ r_1(s, \cdot, z)-r_2(s, \cdot, z))\nonumber\\
&&\quad\quad\quad -F_n(u_1(s-)-u_2(s-))\bigg\}\tilde{N}(ds,dz)\nonumber\\
&&-\int_t^T\int_\RB\bigg \{F_n(u_1(s)-u_2(s)+ r_1(s, \cdot, z)-r_2(s, \cdot, z))-F_n(u_1(s)-u_2(s))\nonumber\\
&&\quad\quad \quad -F_n^{\prime}(u_1(s)-u_2(s))(r_1(s, \cdot, z)-r_2(s, \cdot, z))\bigg\}ds\nu(dz)\nonumber\\
&=:&I^1_n+I^2_n+I^3_n+I^4_n+I^5_n+I^6_n+I^7_n+I^8_n.
\end{eqnarray}
In view of the assumptions (A.1)--(A.3), we have
\begin{eqnarray}\label{2.22}
&&I^2_n=\int_t^TF_n^{\prime}(u_1(s)-u_2(s))(A(s)(u_1(s)-u_2(s)))ds\nonumber\\
&=&\int_t^T\int_Df_n^{\prime}(u_1(s,x)-u_2(s,x))(\sum_{i,j=1}^d {\partial\over {\partial x_i}}(a_{ij}(t,x){\partial\over {\partial x_j}}(u_1(s,x)-u_2(s,x))))dxds  \nonumber\\
&=&-\int_t^T\int_Df_n^{\prime\prime}(u_1(s,x)-u_2(s,x))\sum_{i,j=1}^d a_{ij}(t,x)\frac{\partial}{\partial x_i}(u_1(s,x)-u_2(s,x))\nonumber\\
&&\quad\quad \times  \frac{\partial}{\partial x_j}(u_1(s,x)-u_2(s,x)) dxds,\nonumber\\
&\leq& -\int_t^T (\frac{1}{2a} |\beta|^2(s)+\delta_1)ds\int_Df_n^{\prime\prime}(u_1(s,x)-u_2(s,x))|\nabla (u_1(s,x)-u_2(s,x))|^2dx,\nonumber\\
&&
\end{eqnarray}
\begin{eqnarray}\label{2.23}
I^6_n&=&-\frac{1}{2}\int_t^T\int_Df_n^{\prime\prime}(u_1(s,x)-u_2(s,x))|Z_1(s,x)-Z_2(s,x)|^2dxds,\nonumber\\
&&
\end{eqnarray}
and
\begin{eqnarray}\label{2.001}
&&I^8_n\nonumber\\
&=&-\int_t^T\int_\RB\int_D\bigg \{f_n(u_1(s,x)-u_2(s,x)+ r_1(s, x, z)-r_2(s, x, z))-f_n(u_1(s,x)-u_2(s,x))\nonumber\\
&&\quad\quad \quad -f_n^{\prime}(u_1(s,x)-u_2(s,x))(r_1(s, x, z)-r_2(s, x, z))\bigg\}dxds\nu(dz)\nonumber\\
&=&-\frac{1}{2}\int_t^T\int_\RB \int_D f_n^{\prime\prime}(u_1(s,x)-u_2(s,x)+\theta(s,x,z)(r_1(s, x, z)-r_2(s, x, z)))\nonumber\\
&&\quad\quad\quad \times (r_1(s, x, z)-r_2(s, x, z))^2dxds\nu(dz),
\end{eqnarray}
where $0\leq \theta(s,x,z)\leq 1$.
\vskip 0.3cm
We further write  $ I^3_n$ as
\begin{eqnarray}\label{2.24}
&&I^3_n\nonumber\\
&&=\int_t^T\int_Df_n^{\prime}(u_1(s,x)-u_2(s,x))(b_1(s,x,u_1(s,x),\nabla u_1(s,x), Z_1(s,x), r_1(s,x, \cdot))\nonumber\\
&&\quad\quad\quad -b_2(s,x,u_2(s,x),\nabla u_2(s,x), Z_2(s,x),r_2(s,x, \cdot)))dxds\nonumber\\
&&= \int_t^T\int_Df_n^{\prime}(u_1(s,x)-u_2(s,x))(b_1(s,x,u_1(s,x), \nabla u_1(s,x),Z_1(s,x),r_1(s,x, \cdot))\nonumber\\
&&\quad\quad\quad
-b_2(s,x,u_1(s,x), \nabla u_1(s,x),Z_1(s,x),r_1(s,x, \cdot)))dxds
 \nonumber\\
&&+\int_t^T\int_Df_n^{\prime}(u_1(s,x)-u_2(s,x))(b_2(s,x,u_1(s,x), \nabla u_1(s,x),Z_1(s,x), r_1(s,x, \cdot))\nonumber\\
&&\quad\quad\quad -b_2(s,x,u_2(s,x), \nabla u_1(s,x),Z_1(s,x), r_1(s,x, \cdot)))dxds\nonumber\\
&&+\int_t^T\int_Df_n^{\prime}(u_1(s,x)-u_2(s,x))(b_2(s,x,u_2(s,x), \nabla u_1(s,x),Z_1(s,x), r_1(s,x, \cdot))\nonumber\\
&&\quad\quad\quad -b_2(s,x,u_2(s,x), \nabla u_2(s,x),Z_1(s,x), r_1(s,x, \cdot)))dxds\nonumber\\
&&+\int_t^T\int_Df_n^{\prime}(u_1(s,x)-u_2(s,x))(b_2(s,x,u_2(s,x),\nabla u_2(s,x), Z_1(s,x),r_1(s,x, \cdot))\nonumber\\
&&\quad\quad\quad -b_2(s,x,u_2(s,x), \nabla u_2(s,x), Z_2(s,x),r_1(s,x, \cdot)))dxds\nonumber\\
&&+\int_t^T\int_Df_n^{\prime}(u_1(s,x)-u_2(s,x))(b_2(s,x,u_2(s,x), \nabla u_2(s,x), Z_2(s,x),r_1(s,x, \cdot))\nonumber\\
&&\quad\quad\quad -b_2(s,x,u_2(s,x), \nabla u_2(s,x), Z_2(s,x),r_2(s,x, \cdot)))dxds\nonumber\\
&&:= I^3_{n,1}+I^3_{n,2}+I^3_{n,3}+I^3_{n,4}+I^3_{n,5}.
\end{eqnarray}
Now, $I^3_{n,1}\leq 0$ by the assumption on $b_1$ and $b_2$, and
\begin{equation}\label{2.002}
I^3_{n,2}\leq C\int_t^T\int_D((u_1(s,x)-u_2(s,x))^{+})^2dxds,
\end{equation}
by the Lipschtiz condition of $b_2$. Recalling  the constant $\delta_1$ in (A.1), we can find a constant $C_{\delta_1}$ such that
\begin{eqnarray}\label{2.004}
I^3_{n,3}&\leq &C\int_t^Tds\int_Ddxf_n^{\prime}(u_1(s,x)-u_2(s,x))|\nabla u_1(s,x)-\nabla u_2(s,x)|\nonumber\\
&\leq &\delta_1 \int_t^Tds\int_D f_n^{\prime\prime}(u_1(s,x)-u_2(s,x))|\nabla u_1(s,x)-\nabla u_2(s,x)|^2dx\nonumber\\
&&+C_{\delta_1}\int_t^Tds\int_D \frac{f_n^{\prime}(u_1(s,x)-u_2(s,x))^2}{   f_n^{\prime\prime}(u_1(s,x)-u_2(s,x))}dx\nonumber\\
&\leq &\delta_1 \int_t^Tds\int_D f_n^{\prime\prime}(u_1(s,x)-u_2(s,x))|\nabla u_1(s,x)-\nabla u_2(s,x)|^2dx\nonumber\\
&&+C\int_t^T\int_D ((u_1(s,x)-u_2(s,x))^+)^2dxds,
\end{eqnarray}
where we have used the fact that there exists a constant $C$ (independent of $n$)  such that for $y\geq 0$,
\begin{equation}\label{2.008}
\frac{f_n^{\prime}(x)^2}{f_n^{\prime\prime}(x+y)}\leq C(x^+)^2.
\end{equation}
This can be easily checked using the definition of $f_n$.
By a similar trick, for any $\delta_2>0$, we have
\begin{eqnarray}\label{2.25}
I^3_{n,4}&\leq &C\int_t^T\int_Df_n^{\prime}(u_1(s,x)-u_2(s,x))|Z_1(s,x)-Z_2(s,x)|dxds\nonumber\\
&\leq &\delta_2 \int_t^Tds\int_D f_n^{\prime\prime}(u_1(s,x)-u_2(s,x))|Z_1(s,x)-Z_2(s,x)|^2dx\nonumber\\
&&+C_{\delta_2}\int_t^T\int_D ((u_1(s,x)-u_2(s,x))^+)^2dxds.\nonumber\\
&&
\end{eqnarray}
In view of the assumption (\ref{a3}), we have
\begin{eqnarray}\label{2.003}
I^3_{n,5}&\leq &\int_t^Tds\int_Ddxf_n^{\prime}(u_1(s,x)-u_2(s,x))\int_{\RB}(r_1(s,x, z)-r_2(s,x, z))\lambda(s,z)\nu(dz)\nonumber\\
&\leq& \int_t^Tds\int_Ddxf_n^{\prime}(u_1(s,x)-u_2(s,x))\nonumber\\
&&\quad\quad\times \int_{\RB}(r_1(s,x, z)-r_1(s,x, z))\chi_{\{r_1(s,x,z)>r_2(s,x,z)\}}\lambda(s,z)\nu(dz)\nonumber\\
&\leq &\frac{1}{2}\int_t^T\int_\RB \int_D f_n^{\prime\prime}(u_1(s,x)-u_2(s,x)+\theta(s,x,z)(r_1(s, x, z)-r_2(s, x, z)))\nonumber\\
&&\quad\quad\quad \times (r_1(s, x, z)-r_2(s, x, z))^2\chi_{\{r_1(s,x,z)>r_2(s,x,z)\}}dxds\nu(dz)\nonumber\\
&&+C\int_t^T\int_\RB \lambda(s,z)^2\int_D \frac{f_n^{\prime}(u_1(s,x)-u_2(s,x))^2}{   f_n^{\prime\prime}(u_1(s,x)-u_2(s,x)+\theta(s,x,z)(r_1(s, x, z)-r_2(s, x, z)))}\nonumber\\
&&\quad\quad\quad \quad \times \chi_{\{r_1(s,x,z)>r_2(s,x,z)\}}dxds\nu(dz)\nonumber\\
&\leq &\frac{1}{2}\int_t^T\int_\RB \int_D f_n^{\prime\prime}(u_1(s,x)-u_2(s,x)+\theta(s,x,z)(r_1(s, x, z)-r_2(s, x, z)))\nonumber\\
&&\quad\quad\quad\quad \times (r_1(s, x, z)-r_2(s, x, z))^2dxds\nu(dz)\nonumber\\
&&+C\int_\RB (1\wedge |z|)^2\nu(dz)\int_t^T\int_D ((u_1(s,x)-u_2(s,x))^+)^2dxds,
\end{eqnarray}
where we have used (\ref{2.008}) again.
\vskip 0.3cm
Now, by integration by parts, for  $0<a<1$,
\begin{eqnarray}\label{2.005}
I^4_{n}&=&\int_t^T\int_D<f_n^{\prime}(u_1(s,x)-u_2(s,x))\beta(s), \nabla (Z_1(s,x)-Z_2(s,x))>dxds\nonumber\\
&=&-\int_t^T\int_D<f_n^{\prime\prime}(u_1(s,x)-u_2(s,x))\nabla (u_1(s,x)-u_2(s,x)), \nonumber\\
&&\quad\quad\quad \beta(s)(Z_1(s,x)-Z_2(s,x))>dxds\nonumber\\
&\leq& \frac{1}{2}a\int_t^T\int_Df_n^{\prime\prime}(u_1(s,x)-u_2(s,x))|Z_1(s,x)-Z_2(s,x)|^2dxds\nonumber\\
&&+\frac{1}{2a}\int_t^T\beta^2(s)\int_Df_n^{\prime\prime}(u_1(s,x)-u_2(s,x))|\nabla (u_1(s,x)-u_2(s,x))|^2dxds\nonumber
\end{eqnarray}
Now,  choose $\delta_2>0$ sufficiently  small to satisfy that

\begin{equation}\label{2.007}
\delta_2+\frac{1}{2}a \leq \frac{1}{2}
\end{equation}

Adding (\ref{2.22}),(\ref{2.001}), (\ref{2.004}), (\ref{2.005}), (\ref{2.24}), (\ref{2.25}) and (\ref{2.003}) together  and taking into account of  (\ref{2.007}) we deduce that
\begin{equation}\label{2.26}
I^2_n+I^3_{n}+I^4_{n}+I^6_{n}+I^8_{n}\leq C\int_t^T\int_D((u_1(s,x)-u_2(s,x))^{+})^2dxds
\end{equation}
Thus it follows from (\ref{2.21}), (\ref{2.22}) and (\ref{2.26}) that
 \begin{eqnarray}\label{2.27}
&&F_n(u_1(t)-u_2(t))\nonumber\\
&\leq &F_n(\phi_1-\phi_2)+C\int_t^T\int_D((u_1(s,x)-u_2(s,x))^{+})^2dxds
\nonumber\\
&&-\int_t^TF_n^{\prime}(u_1(s)-u_2(s))(Z_1(s)-Z_2(s))dB_s\nonumber\\
&&-\int_t^T\int_\RB\bigg \{ F_n(u_1(s-)-u_2(s-)+ r_1(s, \cdot, z)-r_2(s, \cdot, z))\nonumber\\
&&\quad\quad -F_n(u_1(s-)-u_2(s-))\bigg\}\tilde{N}(ds,dz)
\end{eqnarray}
Take expectation and let $n\rightarrow \infty$ to get
\begin{equation}\label{2.28}
E[\int_D((u_1(t,x)-u_2(t,x))^{+})^2dx]\leq \int_t^TdsE[\int_D((u_1(s,x)-u_2(s,x))^{+})^2dx]
\end{equation}

Gronwall's  inequality yields that
 \begin{equation}\label{2.29}
E[\int_D((u_1(t,x)-u_2(t,x))^{+})^2dx]=0,
\end{equation}
which completes the proof of the theorem. \quad $\blacksquare$

\section{Application}
Let $u(t,x)$ be the solution of a BSDE of the form (\ref{1.1}) satisfying the conditions (A.1)-(A.3). Assume that $b$ does not
depend on $u$, i.e.
\begin{equation}\label{3.1}
b(t,x,u,\nabla u, Z,r)=b(t,x,Z,r)\quad \mbox{for all}\quad t,x,u,Z,r.
\end{equation}
Moreover, assume that $b(t,x,Z,r)$ is concave with respect to $Z, r$ for all $t, x$. If we, for example, regard $\phi(x)$ as a financial standing at time $t=T$ and at the point $x$, we may as in \cite{OSZ} define the risk $\rho(\phi)(x)$ of $\phi$ at time $t=0$ and at the point $x$ by
\begin{equation}\label{3.2}
\rho(\phi)(x)=-u(0, x); \quad x\in \RB^d.
\end{equation}
Using the comparison theorem (Theorem 2.1) we can now verify that $\phi\rightarrow \rho(\phi)$ is a convex risk measure, in the sense that it satisfies the following conditions:

\noindent (3.3) (Convexity)  $\rho(\lambda\phi_1+(1-\lambda)\phi_2)\leq \lambda\rho(\phi_1)+(1-\lambda)\rho(\phi_2)$ for all $\lambda\in [0,1]$ and
all $\phi_1, \phi_2$.

\noindent (3.4) (Monotonicity)  $\phi_1\leq \phi_2\Rightarrow \rho(\phi_1)\geq \rho(\phi_2)$.

\noindent (3.5) (Translation invariance)  $\rho(\phi+a)=\rho(\phi)-a$ for all $\phi$ and
all constants $a$.
\vskip 0.4cm
Thus we have an extension of the convex risk measure concept (see e.g. \cite{FS}) to a space-dependent situation. This might be of relevance in large systems of interacting components.

\end{document}